\documentclass[letter,10pt]{amsart}
\usepackage{amsmath, amsfonts, amssymb,amscd,graphics,graphicx,color,eucal,setspace} 
\usepackage{latexsym}   
\usepackage{color}

\usepackage{url}

\usepackage{fullpage}

\numberwithin{equation}{section}

\theoremstyle{plain}
\newtheorem{theo}{Theorem}[section]
 
\newtheorem{coro}[theo]{Corollary}

\theoremstyle{definition}
\newtheorem{rema}{Remark}
\newtheorem{exam}{Example}

\def\C{\mathbb C}
\def\Z{\mathbb Z}
\def\Q{\mathbb Q}

\def\x{x}
\def\varphih{\varphi_h}
\def\Ih{I(h)}

\DeclareMathOperator{\Hess}{Hess}

\newcommand{\Flags}{\mathcal{F}\ell ags}

\newcommand{\hsm}{{\hspace{1mm}}}

\begin{document}
  
\title[Cohomology rings of Hessenberg varieties]{The equivariant cohomology rings of regular nilpotent Hessenberg varieties in Lie type A : Research Announcement}

\author {Hiraku Abe}
\address{Osaka City University Advanced Mathematical Institute, Sumiyoshi-ku, Osaka 558-8585, Japan}
\email{hirakuabe@globe.ocn.ne.jp}

\author {Megumi Harada} 
\address{Department of Mathematics and Statistics, McMaster University, 1280 Main Street West, Hamilton, Ontario L8S4K1, Canada}
\email{Megumi.Harada@math.mcmaster.ca}
\urladdr{\url{http://www.math.mcmaster.ca/~haradam}}

\author {Tatsuya Horiguchi}
\address{Department of Mathematics, Osaka City University, Sumiyoshi-ku, Osaka 558-8585, Japan}
\email{d13saR0z06@ex.media.osaka-cu.ac.jp}

\author {Mikiya Masuda}
\address{Department of Mathematics, Osaka City University, Sumiyoshi-ku, Osaka 558-8585, Japan}
\email{masuda@sci.osaka-cu.ac.jp}
\date{\today}



\maketitle

\begin{center} 
\emph{Dedicated to the memory of Samuel Gitler (1933-2014)} 
\end{center}

\begin{abstract}
Let $n$ be a fixed positive integer and $h: \{1,2,\ldots,n\} \rightarrow \{1,2,\ldots,n\}$ a
Hessenberg function. The main result of this manuscript is to give a systematic method for producing 
an explicit presentation by generators and relations of the
equivariant and ordinary cohomology rings (with $\Q$ coefficients)
 of 
any regular nilpotent Hessenberg variety $\Hess(h)$ in type A. Specifically, 
we give an explicit algorithm, depending only on the Hessenberg
function $h$, which produces 
the $n$ defining relations $\{f_{h(j),j}\}_{j=1}^n$ in the
equivariant cohomology ring. Our result generalizes known results: 
for the
case $h=(2,3,4,\ldots,n,n)$, which corresponds to the Peterson
variety $Pet_n$, we recover the presentation of $H^*_S(Pet_n)$ given
previously by Fukukawa, Harada, and Masuda. Moreover, in the case $h=(n,n,\ldots,n)$, for which
the corresponding regular nilpotent Hessenberg variety is the full
flag variety $\Flags(\C^n)$, we can explicitly relate
the generators of our ideal with those in the usual Borel
presentation
of the cohomology ring of $\Flags(\C^n)$. The proof of our main theorem
includes an argument that the restriction homomorphism
$H^*_T(\Flags(\C^n)) \to H^*_S(\Hess(h))$ is surjective.  In this
research announcement, we briefly recount the context and state our
results; we also give a sketch of our proofs and conclude with a
brief discussion of open questions. A manuscript containing more details and full proofs is forthcoming. 
\end{abstract}

\setcounter{tocdepth}{1}
\tableofcontents

\section*{Introduction} \label{sec:intro}

This paper is a research announcement and is a contribution to
the volume dedicated to the illustrious career of Samuel Gitler. A
manuscript containing full details is in preparation
\cite{AbeHaradaHoriguchiMasuda}. 

Hessenberg varieties (in type A) are subvarieties of the full flag
variety $\Flags(\C^n)$ of nested sequences of subspaces in
$\C^n$. Their geometry and (equivariant) topology have been studied extensively since the late 1980s \cite{DeM, DeMPS, DeMShay}. 
This subject 
lies at the intersection of, and makes
connections between, many research areas such as: 
geometric representation theory \cite{Springer76, Fung03},
combinatorics \cite{Fulm, mb}, and algebraic
geometry and topology \cite{br-ca04, InskoYong}. Hessenberg
varieties also arise in the study of the quantum cohomology of the
flag variety \cite{Kostant, Rietsch}. 

The (equivariant) cohomology rings of Hessenberg varieties has been
actively studied in recent years. For instance, Brion and Carrell
showed an isomorphism between the equivariant cohomology ring of a
regular nilpotent Hessenberg variety with the affine coordinate ring
of a certain affine curve \cite{br-ca04}. In the special case of
Peterson varieties $Pet_n$ (in type A), the second author and Tymoczko provided an
explicit set of generators for $H^*_S(Pet_n)$ and also proved a 
Schubert-calculus-type ``Monk formula'', thus giving a presentation of
$H^*_S(Pet_n)$ via generators and relations \cite{HaradaTymoczko}. Using this Monk formula,
Bayegan and the second author derived a ``Giambelli formula'' \cite{BayeganHarada} for
$H^*_S(Pet_n)$ which then yields a simplification of the original presentation
given in \cite{HaradaTymoczko}. Drellich has generalized 
the results in \cite{HaradaTymoczko} and \cite{BayeganHarada} to
Peterson varieties in all Lie types \cite{Drellich}. In another direction, descriptions of the equivariant
cohomology rings of Springer
varieties and regular nilpotent Hessenberg varieties in type A have
been studied by Dewitt and the second author \cite{DewittHarada}, 
the third author \cite{Horiguchi}, the first and third authors \cite{AbeHoriguchi}, and
Bayegan and the second author \cite{BayeganHarada2}.
However, it has been an open question to give a general and systematic
description of the equivariant cohomology rings of all regular
nilpotent Hessenberg varieties \cite[Introduction, page
2]{InskoTymoczko}, to which our results provide an answer (in Lie type
A). 

Finally, we mention that, as a stepping stone to our main result, we can
additionally prove a fact (cf. Section~\ref{sec:sketch}) which seems to be
well-known by experts but for which we did not find an explicit proof
in the literature: namely, that the natural restriction homomorphism
$H^*_T(\Flags(\C^n)) \to H^*_S(\Hess(h))$ is surjective when
$\Hess(h)$ is a regular nilpotent Hessenberg variety (of type A).

\section{Background on Hessenberg varieties} \label{sec:H}

In this section we briefly recall the terminology required to
understand the statements of our main results; in particular we
recall the definition of a regular nilpotent Hessenberg variety, denoted $\Hess(h)$, along with a natural $S^1$-action on it. In this manuscript we only discuss the Lie type A case (i.e. the $GL(n,\C)$ case). We also record some observations regarding the $S^1$-fixed points of $\Hess(h)$, which will be important in later sections.

By the \textbf{flag variety} we mean the homogeneous space
$GL(n,\C)/B$ which may also be identified with 
\[
\mathcal{F}\ell ags(\C^n) := 
\{ V_{\bullet} = (\{0\} \subseteq V_1 \subseteq V_2 \subseteq \cdots V_{n-1} \subseteq
V_n = \C^n) \hsm \mid \hsm \dim_{\C}(V_i) = i \}. 
\]
A \textbf{Hessenberg function} is a function
$h: \{1,2,\ldots, n\} \to \{1,2,\ldots,
n\}$ satisfying $h(i) \geq i$ for all $1 \leq i \leq n$ and $h(i+1)
\geq h(i)$ for all $1 \leq i < n$. We frequently denote a Hessenberg
function by listing its values in sequence, $h = (h(1), h(2), \ldots,
h(n)=n)$. 
Let $N: \C^n \to \C^n$ be a linear operator. 
The \textbf{Hessenberg variety (associated to $N$ and $h$)} 
$\Hess(N,h)$ is defined as the following subvariety of 
$\mathcal{F}\ell ags(\C^n)$: 
\begin{equation}\label{eq:def-Hess}
\Hess(N,h) := \{ V_{\bullet}  \in \mathcal{F}\ell ags(\C^n) \;
\vert \;  N V_i \subseteq
V_{h(i)} \text{ for all } i=1,\ldots,n\} \subseteq \mathcal{F}\ell
ags(\C^n). 
\end{equation}
If $N$ is nilpotent, we say $\Hess(N,h)$ is a \textbf{nilpotent
  Hessenberg variety}, and if $N$ is a principal nilpotent operator
then $\Hess(N,h)$ is
called a \textbf{regular nilpotent Hessenberg variety}. 
In this manuscript we restrict to the regular
nilpotent case, and as such we denote $\Hess(N,h)$ simply as
$\Hess(h)$ where $N$ is understood to be the standard principal
nilpotent operator, i.e. $N$ has one Jordan block with eigenvalue $0$.

Next recall that the following standard torus 

\begin{center}
\begin{equation} \label{eq:T}
T=\left\{ \begin{pmatrix}
  g_1  &    &    &     \\ 
    &  g_2  &    &         \\
    &    &  \ddots  &         \\
    &    &    &      g_n  
\end{pmatrix} \mid  \; g_i\in\mathbb{C}^* \ (i=1,2,\dots n) \right\}
\end{equation} 
\end{center}
acts on the flag variety $Flags(\mathbb{C}^n)$ by left multiplication. However, this $T$-action does not preserve the subvariety $\Hess(h)$ in general. This problem can be rectified by considering instead the action of 
the following circle subgroup $S$ of $T$, which does preserve $\Hess(h)$ (\cite[Lemma 5.1]{ha-ty}): 
\begin{center}
\begin{equation} \label{eq:T}
S :=\left\{ \begin{pmatrix}
  g  &    &    &     \\ 
    &  g^2  &    &         \\
    &    &  \ddots  &         \\
    &    &    &      g^n  
\end{pmatrix} \mid  \; g\in\mathbb{C}^* \right\}. 
\end{equation} 
\end{center}
(Indeed it can be checked that $S^{-1}NS = gN$ which implies that $S$ preserves $\Hess(h)$.) 
Recall that the $T$-fixed points $Flags(\mathbb{C}^n)^T$ of the flag variety $Flags(\mathbb{C}^n)$ 
can be identified with the permutation group $S_n$ on $n$ letters. More concretely, it is straightforward to see that the $T$-fixed points are the set 
$$\{(\langle e_{w(1)}\rangle \subset \langle e_{w(1)},e_{w(2)}\rangle \subset \dots \subset \langle e_{w(1)},e_{w(2)},\dots ,e_{w(n)}\rangle =\mathbb{C}^n) \mid w\in S_n\}$$
where $e_1,e_2,\dots ,e_n$ denote the standard basis of $\mathbb{C}^n$.

It is known that for a regular nilpotent Hessenberg variety $\Hess(h)$ we have 
\[
\mathcal \Hess(h)^S=\mathcal \Hess(h)\cap (Flags(\mathbb{C}^n))^T
\]
so we may view $\Hess(h)^S$ as a subset of $S_n$.

\section{Statement of the main theorem} \label{sec:f}

In this section we state the main result of this paper. 
We first recall some notation and terminology. Let $E_i$ denote the subbundle of the trivial vector bundle $Flags(\mathbb{C}^n)\times \mathbb{C}^n$ over $Flags(\mathbb{C}^n)$ whose fiber at a flag $V_{\bullet}$ is just $V_i$. 
We denote the $T$-equivariant first Chern class of the line bundle $E_i/E_{i-1}$ by $\tilde\tau_i\in H^2_{T}(Flags(\mathbb{C}^n))$. 
Let $\mathbb{C}_i$ denote the one dimensional representation of $T$ through the map $T\rightarrow \mathbb{C}^*$ given by $diag(g_1,\dots,g_n)\mapsto g_i$. 
In addition we denote the first Chern class of the line bundle $ET\times _{T}\mathbb{C}_i$ over $BT$ by $t_i\in H^2(BT)$. 
It is well-known that the $t_1,\dots,t_n$ generate $H^*(BT)$ as a ring and are algebraically independent, so 
we may identify $H^*(BT)$ with the polynomial ring $\mathbb{Q}[t_1,\dots,t_n]$ as rings. 
Furthermore, it is known that $H^{\ast}_{T}(Flags(\mathbb{C}^n))$ is generated as a ring by the elements $\tilde\tau_1,\dots ,\tilde\tau_n,t_1,\dots ,t_n$. Indeed, by sending $x_i$ to $\tilde\tau_i$ and the $t_i$ to $t_i$ we 
obtain the following isomorphism: 
\begin{equation*}\label{eq:equicohofl}
H^{\ast}_{T}(Flags(\mathbb{C}^n))\cong \mathbb{Q}[x_1,\dots,x_n,t_1,\dots,t_n]/(e_i(x_1,\dots ,x_n)-e_i(t_1,\dots ,t_n) \mid 1\leq i\leq n).
\end{equation*}
Here the $e_i$ denote the degree-$i$ elementary symmetric polynomials in the relevant variables. In particular, since the odd cohomology of the flag variety $Flags(\mathbb{C}^n)$ vanishes, we additionally obtain the following: 
\begin{equation}\label{eq:cohofl}
H^{\ast}(Flags(\mathbb{C}^n))\cong \mathbb{Q}[x_1,\dots,x_n]/(e_i(x_1,\dots ,x_n) \mid 1\leq i\leq n).
\end{equation}
As mentioned in Section~\ref{sec:H}, in this manuscript we focus on a particular circle subgroup $S$ of the usual maximal torus $T$. For this subgroup $S$, we denote the first Chern class of the line bundle $ES\times _{S}\mathbb{C}$ over $BS$ by $t\in H^2(BS)$, where by $\mathbb{C}$ we mean the standard one-dimensional representation of $S$ through the map $S\rightarrow \mathbb{C}^*$ given by $diag(g,g^2,\dots,g^n)\mapsto g$. 
Analogous to the identification $H^*(BT) \cong \mathbb{Q}[t_1, \ldots, t_n]$, we may also identify 
$H^*(BS)$ with $\mathbb{Q}[t]$ as rings.

Consider the restricion homomorphism 
\begin{equation}\label{eq:rest}
H^*_T(\Flags(\mathbb{C}^n)) \to H^*_S(\Hess(h)).
\end{equation}
Let $\tau_i$ denote the image of $\tilde\tau_i$ under \eqref{eq:rest}.
We next analyze some algebraic relations satisfied by the $\tau_i$. For this purpose, we now introduce some polynomials  $f_{i,j}=f_{i,j}(x_1,\dots,x_n,t)\in\mathbb{Q}[x_1,\dots,x_n,t]$.

First we define 
\begin{equation} \label{eq:f-1}
p_i:=\sum_{k=1}^i(x_k-kt)\quad (1\leq i\leq n).
\end{equation} 
For convenience we also set $p_0 :=0$ by definition. 
Let $(i,j)$ be a pair of natural numbers satisfying  $n \geq i \geq j \geq 1$. These polynomials should be visualized as being associated to the $(i,j)$-th spot in an $n \times n$ matrix. Note that by assumption on the indices, we only define the $f_{i,j}$ for entries in the lower-triangular part of the matrix, i.e. the part at or below the diagonal. The definition of the $f_{i, j}$ is inductive, beginning with the case when $i=j$, i.e. the two indices are equal. 
In this case we make the following definition: 
\begin{equation} \label{eq:f-2}
f_{j,j}:=p_j \quad (1\leq j\leq n). 
\end{equation}
Now we proceed inductively for the rest of the $f_{i,j}$ as follows: for $(i,j)$ with $n\ge i>j\ge 1$ we define: 
\begin{equation} \label{eq:f-3}
f_{i,j}:=f_{i-1,j-1}+\big(x_j-x_i-t\big)f_{i-1,j}.
\end{equation}
Again for convenience we define $f_{*,0}:=0$ for any $*$. 
Informally, we may visualize each $f_{i,j}$ as being associated to the lower-triangular $(i,j)$-th entry in an $n \times n$ matrix, as follows: 
\begin{equation}\label{eq:matrix fij}
\begin{pmatrix}
f_{1,1} & 0 & \cdots & \cdots & 0 \\
f_{2,1} & f_{2,2} & 0 & \cdots &  \\
f_{3,1} & f_{3,2} & f_{3,3} & \ddots & \\
\vdots &        &                 &  &  \\
f_{n,1} & f_{n,2} & \cdots &  & f_{n,n} 
\end{pmatrix}
\end{equation}

To make the discussion more concrete, we present an explicit example. 

\begin{exam}\label{exam:f}
Suppose $n=4$. Then the $f_{i,j}$ have the following form. 

\noindent
$f_{i,i}=p_i$ \ ($1\leq i\leq 4$)

\noindent
$f_{2,1}=(x_1-x_2-t)p_1$

\noindent
$f_{3,2}=(x_1-x_2-t)p_1+(x_2-x_3-t)p_2$

\noindent
$f_{4,3}=(x_1-x_2-t)p_1+(x_2-x_3-t)p_2+(x_3-x_4-t)p_3$

\noindent
$f_{3,1}=(x_1-x_3-t)(x_1-x_2-t)p_1$

\noindent
$f_{4,2}=(x_1-x_3-t)(x_1-x_2-t)p_1+(x_2-x_4-t)\{(x_1-x_2-t)p_1+(x_2-x_3-t)p_2\}$

\noindent
$f_{4,1}=(x_1-x_4-t)(x_1-x_3-t)(x_1-x_2-t)p_1$

\medskip
\end{exam}

For general $n$, the polynomials $f_{i,j}$ for each $(i,j)$-th entry in the matrix~\eqref{eq:matrix fij} above can also be expressed 
in a closed formula in terms of 
certain polynomials $\Delta_{i, j}$ for $i\geq j$ which are determined inductively, starting on the main diagonal.  
As for the $f_{i,j}$, we think of $\Delta_{i,j}$ for $i \geq j$ as being 
associated to the $(i,j)$-th box in an $n \times n$ matrix. 
In what follows, for $0<k\leq n-1$, we refer to the lower-triangular matrix entries in the $(i,j)$-th spots where $i-j=k$ as the \textbf{$k$-th lower diagonal}. (Equivalently, the $k$-th lower diagonal is the ``usual'' diagonal of the lower-left $(n-k) \times (n-k)$ submatrix.) The usual diagonal is the $0$-th lower
diagonal in this terminology. We now define the $\Delta_{i,j}$ as follows. 
\begin{enumerate} 
\item First place the linear polynomial $x_i - it$ in the $i$-th entry along the $0$-th lower (i.e. main) diagonal, so $\Delta_{i,i} := x_i - it$. 
\item Suppose that $\Delta_{i,j}$ for the $k-1$-st lower diagonal have already been defined. 
Let $(i,j)$ be on the $k$-th lower diagonal, so $i-j=k$. Define 
\begin{equation*}
\Delta_{i,j} := \left(\sum_{\ell=1}^j \Delta_{i-j+\ell-1, \ell}\right)(x_j - x_i - t).
\end{equation*}
\end{enumerate} 
In words, this means the following. Suppose $k=i-j>0$. Then $\Delta_{i,j}$ is the product of $(x_j - x_i -t)$ with the sum of the entries in the boxes 
which are in the ``diagonal immediately above the $(i,j)$ box'' (i.e. the boxes which are in the $(k-1)$-st lower diagonal), but we omit any boxes to the right of the $(i,j)$ box (i.e. in columns $j+1$ or higher).  
Finally, the polynomial $f_{i,j}$ is obtained by 
taking the sum of the entries in the $(i,j)$-th box and any boxes ``to its left'' in the same lower diagonal. More precisely, 
\begin{equation} \label{eq:Delta-4}
f_{i,j}=\sum_{k=1}^j\Delta_{i-j+k,k}.
\end{equation}

We are now ready to state our main result. 

\begin{theo} \label{theo:main}
Let $n$ be a positive integer and $h: \{1,2,\ldots,n\} \to \{1,2,\ldots,n\}$ a Hessenberg function. Let $\Hess(h) \subset \Flags(\mathbb{C}^n)$ denote the corresponding regular nilpotent Hessenberg variety equipped with the circle $S$-action described above. Then the restriction map 
\begin{equation*}
H^*_T(\Flags(\mathbb{C}^n)) \to H^*_S(\Hess(h))
\end{equation*} 
is surjective.
Moreover, there is an isomorphism of $\mathbb{Q}[t]$-algebras
\begin{equation*}
H^{\ast}_{S}(\mathcal \Hess(h))\cong \mathbb{Q}[x_1,\dots,x_n,t]/\Ih 
\end{equation*}
sending $x_i$ to $\tau_i$ and $t$ to $t$ and we identify $H^{\ast}(BS)=\mathbb{Q}[t]$. Here the ideal $\Ih$ is defined by 
\begin{equation}\label{eq:def Ih} 
\Ih:=(f_{h(j),j} \mid 1\leq j\leq n).
\end{equation} 
\end{theo}

We can also describe the ideal $\Ih$ defined in~\eqref{eq:def Ih} as follows. 
Any Hessenberg function $h: \{1,2,\ldots,n\} \to \{1,2,\ldots,n\}$ determines a subspace of the vector space $M(n \times n,\mathbb{C})$ of matrices as follows: an $(i,j)$-th entry is required to be $0$ if $i>h(j)$. If we represent a Hessenberg function $h$ by listing its values $(h(1), h(2), \cdots, h(n))$, then the Hessenberg subspace can be described in words as follows: the first column (starting from the left) is allowed $h(1)$ non-zero entries (starting from the top), the second column is allowed $h(2)$ non-zero entries, et cetera. For example, if $h = (3,3,4,5,7,7,7)$ then the Hessenberg subspace is 
\[
\left\{ \begin{pmatrix} 
\star & \star & \star & \star & \star & \star& \star \\
\star & \star & \star & \star & \star & \star& \star \\
\star & \star & \star & \star & \star & \star& \star \\
0 & 0 &  \star & \star & \star & \star& \star \\
0 & 0 &  0 & \star & \star & \star& \star \\
0 & 0 &  0 & 0 & \star & \star& \star \\
0 & 0 &  0 & 0 & \star & \star& \star \\
\end{pmatrix} 
\right\} \subseteq M(7 \times 7,\mathbb{C}).
\]
Then, using the association of the polynomials $f_{i,j}$ with the $(i,j)$-th entry of the matrix~\eqref{eq:matrix fij}, the ideal $\Ih$ can be described as being ``generated by the $f_{i,j}$ in the boxes at the bottom of each column in the Hessenberg space''. For instance, in the $h=(3,3,4,5,7,7,7)$ example above, the generators are $\{f_{3,1}, f_{3,2}, f_{4,3}, f_{5,4}, f_{7,5}, f_{7,6}, f_{7,7}\}$. 

Our main result generalizes previous known results. 

\begin{rema}\label{rema:peterson}
Consider the special case $h=(2,3,\dots,n,n)$. In this case the corresponding regular nilpotent Hessenberg variety has been well-studied and it is called a \textbf{Peterson variety} $Pet_n$ (of type $A$). Our result above is a generalization of the result in \cite{fu-ha-ma} which gives a presentation of $H^*_S(Pet_n)$. Indeed, for $1\leq j\leq n-1$, we obtain from \eqref{eq:f-3} and \eqref{eq:f-1}
that 
\begin{equation*}
\begin{split}
f_{j+1,j}&=f_{j,j-1}+(x_j-x_{j+1}-t)f_{j,j} \\
         &=f_{j,j-1}+(-p_{j-1}+2p_j-p_{j+1}-2t)p_j
\end{split}
\end{equation*}
and since $f_{n,n}=p_n$ we have 
\begin{equation*}
\begin{split}
H^{\ast}_{S}(Pet_n)&\cong \mathbb{Q}[x_1,\dots,x_n,t]/\big(f_{j,j-1}+(-p_{j-1}+2p_j-p_{j+1}-2t)p_j, \ p_n \mid 1\leq j\leq n-1\big) \\
                 &= \mathbb{Q}[x_1,\dots,x_n,t]/\big((-p_{j-1}+2p_j-p_{j+1}-2t)p_j, \ p_n \mid 1\leq j\leq n-1\big) \\
                 &\cong \mathbb{Q}[p_1,\dots,p_{n-1},t]/\big((-p_{j-1}+2p_j-p_{j+1}-2t)p_j \mid 1\leq j\leq n-1\big)
\end{split}
\end{equation*}
which agrees with \cite{fu-ha-ma}. (Note that we take by convention $p_0=p_n=0$.) 
\end{rema}

The main theorem above also immediately yields a computation of the ordinary cohomology ring. Indeed, since the 
odd degree cohomology groups of $\Hess(h)$ vanish \cite{ty}, by setting $t=0$ we obtain the ordinary cohomology. Let $\check{f}_{i,j} := f_{i,j}(x,t=0)$ denote the polynomials in the variables $x_i$ obtained by setting $t=0$. A computation then shows that 
\[
\check{f}_{i,j} = \sum_{k=1}^j x_k \prod_{\ell=j+1}^i (x_k - x_\ell).
\]
(For the case $i=j$ we take by convention $\prod_{\ell=j+1}^i (x_k-x_\ell)=1$.) We have the following.

\begin{coro} \label{coro:main}
Let the notation be as above. There is a ring isomorphism
\begin{equation*}
H^{\ast}(\mathcal \Hess(h))\cong \mathbb{Q}[x_1,\dots,x_n]/\check{I}(h)
\end{equation*}
where $\check{I}(h) :=\big(\check{f}_{h(j), j} \mid 1\leq j\leq n\big)$.
\end{coro}

\begin{rema}\label{rema:full flag} 
Consider the special case $h=(n,n,\ldots,n)$. In this case the condition in~\eqref{eq:def-Hess} is vacuous and the associated regular nilpotent Hessenberg variety is the full flag variety $\Flags(\mathbb{C}^n)$. In this case we can explicitly relate the generators $\check{f}_{h(j)=n, j}$ of our ideal $\check{I}(h) = \check{I}(n,n,\ldots,n)$ with 
the power sums $\textsf{p}_r(x) = \textsf{p}_r(x_1, \ldots, x_n) := \sum_{k=1}^n x_k^r$, thus relating our presentation with the usual Borel presentation as in~\eqref{eq:cohofl}, see e.g. \cite{fult97}. More explicitly, 
for $r$ be an integer, $1\le r\le n$, define 
\[
\textsf{q}_r(\x)=\textsf{q}_r(\x_1,\dots,\x_n):=\sum_{k=1}^{n+1-r} \x_k\prod_{\ell=n+2-r}^n(\x_k-\x_\ell).
\]
Note that by definition $\textsf{q}_r(\x) = \check{f}_{n, n+1-r}$ so these are the generators of $\check{I}(n,n,\ldots,n)$. The polynomials $\textsf{q}_r(x)$ and the power sums $\textsf{p}_r(x)$ can then be shown to satisfy the relations 
\begin{equation}\label{eq:q and p}
\textsf{q}_r(\x)=\sum_{i=0}^{r-1}(-1)^ie_i(\x_{n+2-r},\dots,\x_n)\textsf{p}_{r-i}(\x). 
\end{equation}

\end{rema} 

\begin{rema} 
In the usual Borel presentation of $H^*(\Flags(\mathbb{C}^n))$, the ideal $I$ of relations is taken to be generated by the elementary symmetric polynomials.  The power sums $\mathsf{p}_r$ generate this ideal $I$ when we consider the cohomology with $\mathbb{Q}$ coefficients, but this is not true with $\Z$ coefficients. Thus our main Theorem~\ref{theo:main} does not hold with $\Z$ coefficients in the case when $h=(n,n,\ldots,n)$, suggesting that there is some subtlety in the relationship between the choice of coefficients and the choice of 
generators of the ideal $I(h)$. 
\end{rema}

\section{Sketch of the proof of the main theorem}\label{sec:sketch}

We now sketch the outline of the proof of the main result (Theorem~\ref{theo:main}) above. 
As a first step, we show that the elements $\tau_i$ satisfy the relations $f_{h(j),j} = f_{h(j),j}(\tau_1, \ldots, \tau_n, t)=0$. The main technique of this part of the proof is (equivariant) localization, i.e. the injection 
\begin{equation}\label{eq:localization}
H^*_S(\Hess(h)) \to H^*_S(\Hess(h)^S). 
\end{equation}
Specifically, we show that the restriction $f_{h(j),j}(w)$ of each $f_{h(j),j}$ to an $S$-fixed point $w \in \Hess(h)^S$ is equal to $0$; by the injectivity of~\eqref{eq:localization} this then implies that $f_{h(j),j}=0$ as desired. This part of the argument is rather long and requires a technical inductive argument based on a particular choice of total ordering on $\Hess(h)^S$ which refines a certain natural partial order on Hessenberg functions. Once we show $f_{h(j),j}=0$ for all $j$, we obtain a well-defined ring homomorphism which sends $x_i$ to $\tau_i$ and $t$ to $t$: 
\begin{equation}\label{eq:homo}
\varphi_h: \mathbb{Q}[x_1,\dots,x_n,t]/(f_{h(j),j} \mid 1\leq j\leq n)\rightarrow H^{\ast}_{S}(\Hess(h)).
\end{equation}
We then show that the two sides of~\eqref{eq:homo} have identical Hilbert series. This part of the argument is rather 
straightforward, following the techniques used in e.g. \cite{fu-ha-ma}.

The next key step in our proof of Theorem~\ref{theo:main} relies on the following two key ideas: firstly, we use our knowledge of the special case where the Hessenberg function $h$ is $h=(n,n,\ldots,n)$, for which the associated regular nilpotent Hessenberg variety is the full flag variety $\Flags(\mathbb{C}^n)$, and secondly, we 
consider localizations of the rings in question with respect to $R:=\Q[t]\backslash\{0\}$. 
For the following, for $h=(n,n,\ldots,n)$ we let $\mathcal{H} := \Hess(h=(n,n,\ldots,n)) = \Flags(\mathbb{C}^n)$ denote the full flag variety and let $I$ denote the associated ideal $I(n,n,\ldots,n)$. 
In this case we know that the map $\varphi := \varphi_{(n,n,\ldots,n)}$ is surjective since the Chern classes $\tau_i$ are known to generate the cohomology ring of $\Flags(\mathbb{C}^n)$. Since the Hilbert series of both sides are identical, we then know that $\varphi$ is an isomorphism. 

The following commutative diagram is crucial for the remainder of the argument. 
\[
\begin{CD}
R^{-1}\big(\Q[x_1,\dots,x_n,t]/I\big) @>R^{-1}\varphi >\cong> R^{-1}H^*_S(\mathcal H) @>>\cong> R^{-1}H^*_S(\mathcal H^S) \\
@VV\text{surj}V @VVV @VV\text{surj}V \\
 R^{-1}\big(\Q[x_1,\dots,x_n,t]/\Ih\big) @>R^{-1}\varphih >> R^{-1}H^*_S(\Hess(h)) @>>\cong> R^{-1}H^*_S(\Hess(h)^S) 
\end{CD}
\]
The horizontal arrows in the right-hand square are isomorphisms by the localization theorem. Since $\varphi$ is an isomorphism, so is $R^{-1}\varphi$. The rightmost and leftmost vertical arrows are easily seen to be surjective, implying that $R^{-1}\varphi_h$ is also surjective. A comparison of Hilbert series shows that $R^{-1}\varphi_h$ is an isomorphism.  
Finally, to complete the proof we consider the commutative diagram
\[
 \begin{CD}
 \Q[x_1,\dots,x_n,t]/\Ih @>\varphih>> H^*_S(\Hess(h)) \\
 @VV\text{inj}V @VV\text{inj}V \\
 R^{-1}\Q[x_1,\dots,x_n,t]/\Ih @>R^{-1}\varphih >\cong > R^{-1}H^*_S(\Hess(h))
 \end{CD}
 \]
 for which it is straightforward to see that the vertical arrows are injections. From this it follows that $\varphi_h$ is an injection, and once again a comparison of Hilbert series shows that $\varphi_h$ is in fact an isomorphism.

 \section{Open questions} 
 
 We outline a sample of possible directions for future work. 
 
 \begin{itemize} 
 
 \item In \cite{mb-ty13}, Mbirika and Tymoczko suggest a possible presentation of the cohomology rings of regular nilpotent Hessenberg varieties. Using our presentation, we can show that the Mbirika-Tymoczko ring is not isomorphic to $H^*(\Hess(h))$ in the special case of Peterson varieties for $n-1\geq 2$, i.e. when $h(i)=i+1, 1 \leq i < n$ and $n\geq 3$. (However, they do have the same Betti numbers.) 
In the case $n=4$, we have also checked explicitly for the Hessenberg functions $h=(2,4,4,4), h=(3,3,4,4)$, and $h=(3,4,4,4)$ that the relevant rings are not isomorphic. 
It would be of interest to understand the relationship between the two rings in some generality. 
 
 \item In \cite{ha-ho-ma}, the last three authors give a presentation of the (equivariant) cohomology rings of Peterson varieties for general Lie type in a pleasant uniform way, using entries in the Cartan matrix. It would be interesting to give a similar uniform description of the cohomology rings of regular nilpotent Hessenberg varieties for all Lie types. 
 
 \item In the case of the Peterson variety (in type A), a basis for the $S$-equivariant cohomology ring was found by the second author and Tymoczko in \cite{HaradaTymoczko}. In the general regular nilpotent case, and following ideas of the second author and Tymoczko \cite{ha-ty}, it would be of interest to construct similar additive bases for $H^*_S(\Hess(h))$. Additive bases with suitable geometric or combinatorial properties could lead to an interesting `Schubert calculus' on regular nilpotent Hessenberg varieties. 
 

\item Fix a Hessenberg function $h$ and let 
$S: \C^n \to \C^n$ be a \emph{regular semisimple} linear operator, i.e. a diagonalizable operator with distinct eigenvalues. 
There is a natural Weyl group action on 
the cohomology ring $H^*(\Hess(S,h))$ of the regular semisimple Hessenberg variety corresponding to $h$ (cf. for instance 
\cite[p. 381]{Tymoczko-Perm} and also \cite{Teff}). Let $H^*(\Hess(S, h))^W$ denote the ring of $W$-invariants where $W$ denotes the Weyl group.  
It turns out that there exists a surjective ring homomorphism $H^*(\Hess(N, h)) \to H^*(\Hess(S, h))^W$ which is an isomorphism in the special case of the Peterson variety. (Historically this line of thought goes back to Klyachko's 1985 paper \cite{Klyachko}.) In an ongoing project, we are investigating properties of this ring homomorphism for general Hessenberg functions $h$.

 \end{itemize}

\end{document}